\documentclass[12pt,a4paper,twoside, english]{amsart}
\usepackage{amssymb}
\usepackage{amscd}
\usepackage{enumitem}
\usepackage{amsmath}

\usepackage[left=2.5cm,right=3.5cm]{geometry}

\usepackage{mathtools}
\usepackage{amsfonts}
\usepackage{amsthm}
\usepackage{graphicx}
\usepackage{mathtools}
\usepackage[colorlinks=true, allcolors=blue]{hyperref}
\usepackage{setspace}
\usepackage{lipsum}
\usepackage{biblatex}
\usepackage[utf8]{inputenc}

\newtheorem{theorem}{Theorem}[section]

\newtheorem{lemma}[theorem]{Lemma}

\newtheorem{corollary}[theorem]{Corollary}

\theoremstyle{definition}

\newtheorem{definition}[theorem]{Definition}

\newtheorem{claim}[theorem]{Claim}

\newtheorem{remark}[theorem]{Remark}

\newtheorem*{theorem*}{Theorem}
\newtheorem*{lemma*}{Lemma}

\newcommand{\ack}{{\bf Acknowledgements}}
\newcommand{\prf}{{\bf Proof}}
\newcommand{\mF}{\mathbb F}
\newcommand{\mE}{\mathbb E}
\newcommand{\mN}{\mathbb N}
\newcommand{\mC}{\mathbb C}
\newcommand{\mZ}{\mathbb Z}
\newcommand{\mR}{\mathbb R}
\newcommand{\mT}{\mathbb T}

\makeatletter
\newcommand{\dotDelta}{{\vphantom{\Delta}\mathpalette\d@tD@lta\relax}}
\newcommand{\d@tD@lta}[2]{%
  \ooalign{\hidewidth$\m@th#1\mkern-1mu\cdot$\hidewidth\cr$\m@th#1\Delta$\cr}%
}
\makeatother

\title[Orthogonality of the M\"obius function to polynomials]{Orthogonality of the M\"obius function to polynomials with applications to Linear Equations in Primes over $\mF_p[x]$}
\author{Tal Meilin}
\makeatother

\begin{document}
\maketitle
\begin{abstract}
    We prove that the M\"obius function is orthogonal to polynomials over $\mF_q[x]$ (up to a characteristic condition). We use this orthogonality property to count prime solutions to affine-linear equations of bounded complexity in $\mathbb{F}_p[x]$, with analog to the work of Green and Tao in \cite{LEP0}. 
\end{abstract}
\tableofcontents

\section{Introduction}\label{intr}
Let  $\mu:\mF_q[x]\rightarrow\mC$ be the M\"obius function, defined by:
\begin{align*}
 \mu(g)=
\begin{dcases}
(-1)^k  &\text{if $k$ is the number of monic irreducible}\\
        &\text{factors of  a square-free $g$}\\
0       &\text{otherwise}
\end{dcases}   
\end{align*}

Our main result in this paper is the following orthogonality property: 

\begin{theorem}[M\"obius function is orthogonal to Polynomials]\label{main}
For any given $m\in\mN$, there exists some  $0<t_m<1$ and $N_m\in \mN$, such that the following holds: for every prime  $p>2m$, $q:=p^s$, any degree $m$ polynomial form $P:\mF_q[x] \to \mF_q$, any character $\chi$ of  $ \mF_q$ and $n>N_m$, we have:
\[
|\mE_ {g\in G_n}\mu(g)\chi(P(g))| \leq q^{-n^{t_m}}.
\]

\end{theorem}
As a consequence, we obtain the following result about solutions of linear equations in primes in finite field geometry:
\begin{theorem}[Linear Equations in Primes in {$\mF_p[x]$}]\label{LEP}
Let $\Psi:\mF_p[x]^d\rightarrow\mF_p[x]^t$ be an affine linear form with complexity $s$ smaller than  $\frac{p}{2}$ (defined in \ref{s-com}), satisfying that for every $i$
    \[\psi_i((g_1,\ldots,g_d))=\sum_{j=1}^d L_{i,j}g_j+b_i \text{ for some }L_{i,j}, b_j\in\mF_p[x]\]
Then:
\[\sum_{g\in G_n^d}\prod_{i\in[t]}\Lambda(\psi_i(g))=p^{nd}\prod_{r \text{ irreducible monics}} \beta_r+o_{d,t,L,s}(p^{nd})\]
Where $\{\beta_r \}$ are the local factors defined in \ref{local_factor}  , and $G_n$ is the additive group of polynomials of degree smaller than $n$, and $$L=\max\{\deg(L_{i,j}), \deg(b_i)| 1\leq i\leq t, 1\leq j\leq d\}.$$
\end{theorem}
Theorem \ref{main} is a generalization of the work of Bienvenu and Ho{\`a}ng L{\^e}, that proved a very similar theorem for quadratic polynomials in \cite{BHL}, and Porrit, that proved the linear case independently \cite{Porrit}.

Theorem \ref{LEP} is the finite field geometry analogue of the main theorem in \cite{LEP0} (Theorem 1.2). The result in \cite{LEP0} was proved conditionally; it relied on the following two conjectures: the Inverse conjecture for the Gowers Norms and the M\"obius and Nilsequence conjecture. Since then, both conjectures have been proven: the first one by Green, Tao and Ziegler  \cite{IGN1}, and the second one by  Green and Tao \cite{Nil}. These conjectures, along with a transference principle introduced in \cite{LAP}, are the main foundations that were needed to establish the main theorem in \cite{LEP0}.

In the $\mF_q[x]$ context, the Inverse Theorem for the Gowers Norms has already been proven by Tao and Ziegler in \cite{IGN}. The transference principle for function fields was proven by  Ho{\`a}ng L{\^e} in \cite{HLFF}. 

Hence, the remaining challenge in proving Theorem \ref{LEP} is Theorem \ref{main}, which is the finite field geometry analogue to the orthogonality of the (arithmetic) M\"obius  function to nil-sequences  ({\cite{Nil}). We prove Theorem \ref{main} in the first part of our paper.
In the second part, we use Theorem \ref{main} to conclude Theorem \ref{LEP}. \\

Besides Theorem \ref{LEP}, the following are some immediate results that can be obtained by Theorem \ref{main} together with  the inverse theorem for the Gowers norms over finite fields (\cite{IGN}):
\begin{theorem}\label{gn0}
    Let $\mF_p$ be a prime field, and 
    $\mu$ be  the M\"obius function $\mu:\mF_p[x]\rightarrow\{-1,0,1\}$. Then, for every $d<\frac{p}{2}$:
    \[||\mu||_{U^d(G_n)}\rightarrow 0.\]
\end{theorem}
where $U^d(G_n)$ is the $d$th-Gowers unifromity norm defined in \ref{gn}. From this, we may deduce an estimate for the correlations of the M\"obius function along linear forms:
\begin{corollary}\label{lfm}
 For a linear form $l_1,\ldots, l_k:\mF_p[x]^d\rightarrow \mF_p[x]$  of complexity $<\frac{p}{2}$ (defined in \ref{s-com})
 \[|\mE_{g\in (G_n)^d}\prod_{i=1}^k\mu(l_i(g))|=o(1).\]
\end{corollary}
In particular, since the linear-form 
$L=(l_1,\ldots,l_k):\mF_q[x]^2\rightarrow\mF_q[x]$ defined by:
\[l_i(f,g)=f+ig\]
and has complexity $k-1$, we get the following result regarding arithmetic progressions:
\begin{theorem} Fix $k \in \mN$. Then for every prime number $p>2k$ we have:
\[\mE_{f,g\in G_n}\mu(f)\mu(f+g)\cdot \ldots \cdot\mu(f+kg)\xrightarrow{} 0.\]
\end{theorem}
Note that with some stronger conditions on $
q$ ,  Sawin and Shusterman proved a stronger version of \ref{lfm} in \cite{SaSh}:
\begin{theorem}[Sawin-Shusterman \cite{SaSh}]\label{bla}
For an odd prime number $p$,  an integer $k>1$, and a power $q$ of $p$ satisfying $q>p^2k^2e^2$, the following holds:\\
For distinct $h_1,\ldots, h_k\in \mF_q[x]$ we have:
\[\sum_{f\in A_n}\mu(f+h_1)\ldots\mu(f+h_k)=o(q^n), n\xrightarrow{}\infty
\]
\end{theorem}
We remark that while Theorem \ref{bla} is much stronger than Theorem \ref{main}, it heavily relies on the specific properties of the M\"obius function. 

Since the only requirement of Theorems \ref{gn0} and \ref{lfm} is the orthogonality of $\mu$ to polynomials (meaning, ineffective bound would be sufficient) we intend to prove similar results for a much larger set of arithmetic functions in a followup paper.\\

\ack. I would like to thank my academic advisor Prof. Tamar Ziegler for her guidance and support throughout this work. I would also like to thank Dr. Joni Teräväinen for his helpful comments and suggestions. The author is supported by ISF grant 2112/20.

\section{Overview of the paper}
In our paper, we establish two major theorems: Theorem \ref{main} and Theorem \ref{LEP}.\\

In the first sections (\ref{intr}-\ref{auxiliary}), we present our main results, along with notations and previous results we'll use in our proof.

In section \ref{proof of theorem 1}, we prove  Theorem \ref{main}, which is the main missing result for the proof of Theorem \ref{LEP} in the context of function fields. We already know that Theorem \ref{main} holds for linear and quadratic polynomials (\cite{BHL}). We prove the theorem by induction over the polynomial's degree. The induction base is the linear case that was proved in \cite{BHL}, and the proof of the step will be done by reaching a contradiction, namely assuming the existence of a polynomial $P$ of degree $m$ that doesn't satisfy Theorem \ref{main}.  To gain some initial information about $P$, we use a version of Vaughan's Lemma For function fields  (Lemma \ref{VL}). In \cite{BHL}, a similar Lemma was used to deduce the quadratic case from the linear case. We aim to prove that such $P$ must have a relatively low rank (see definition\ref{AR}). Intuitively, that would imply that $P$ is a "simple" combination of a "small" number of polynomials of lower degree, which allows us to utilize the induction's assumption. Hence (by   Lemma \ref{IL}) we can conclude that $P$ satisfies Theorem \ref{main}.

While proving that $P$ has a low rank, we apply equidistribution theorems to the results that were obtained by using Lemma \ref{VL} to polynomials over finite fields (Theorem \ref{BR}). 
For large enough $n$ (depending on the polynomial's degree) we can find an $m$-linear map related to $P$ that has a low partition rank on some subspace of $G_n$.  By repetitively using Theorem \ref{BR} along with some other polynomials' properties appearing in section \ref{auxiliary}, we can deduce that this $m$-linear map can be extended to the restriction of $Q_P$ (the $m$-derivative of $P$) to some subspace of $G_n$ with a small co-dimension. This will imply that the rank of $P$ is low as we intended.

We use the equivalence between the low Schmidt rank of $P$ and the low partition rank of its corresponding $m-linear$ form quite often in our proof. The restriction on the characteristic comes from the fact that this equivalence is only valid if $m<Char(\mF_q)$.\\

In section \ref{proof of lep}, we use Theorem \ref{main} to prove Theorem \ref{LEP}. As mentioned before, the proof is an analogy of the classical case's proof in \cite{LEP0}. 

First, to avoid obstacles caused by the behavior of $\Lambda$ on small irreducible polynomials, we replace $\Lambda$ with a correspondent function $\Lambda_{W,b}$, which "ignores" small irreducible polynomials. This strategy is used in \cite{LEP0} and is called "The W-trick".
Then, by using theorem \ref{main}, we exploit the relations between $\Lambda$ and $\mu$ to prove that $\Lambda_{W,b}-1$ is orthogonal to polynomials with degree at most $\frac{p}{2}$. 
By using some properties of $\Lambda$ as well as Theorem \ref{tp} and Theorem \ref{ign} we may deduce that $\Lambda_{W,b}-1$ has small Gowers Norms (see \ref{gn}), which means  that its average on linear forms vanishes (see \ref{VN}). 
Finally, we reverse our focus once again from $\Lambda_{W,b}-1$ to $\Lambda$. Since we proved that the average of $\Lambda_{W,b}-1$ on linear forms vanishes, we can obtain the desired result for $\Lambda$ with the help of some additional calculations.

\section{Notations}
\label{notations}
We start with some basic definitions and notations regarding polynomials over $\mF_q[x].$

For an $m$-linear map $Q:\mF_q[x]^m \to \mF_q$ (i.e, Q is linear in every variable), let $P_Q$ be the homogeneous polynomial corresponding to $Q$, obtained by restricting $Q$ to the diagonal, i.e \[P_Q(g):=Q(g,\ldots,g),\text{ for every } g \in \mF_q[x].\] 
For a homogeneous polynomial  $
P $ of degree $m$ over $ \mF_q[x]$ , let $Q_P$ be the $m$-linear map corresponding to $P$, i.e.  
\[Q_P(h_1,\ldots,h_m):=\Delta_{h_1},\ldots,\Delta_{h_m}P(g)\]
where $h_i \in \mF_q[x]$ and $\Delta_h P(g):=P(g+h)-P(g)$. Note that $Q_P$ is independent of $g$ and:
\[P_{Q_P}=m! P \text{  and  } Q_{P_Q}=m! Q.\]
We state now some important definitions:
\begin{definition}[Partition Rank]\label{PR}
For an $m$-linear map $Q$ define the partition rank - $r_{pr}(Q)$ of $Q$  as the minimal $r$ s.t:
\[
Q= \sum_{i=1}^r M_iR_i
\]
where for $1 \le i \le r$, $M_i, R_i$ are $k_i$-linear and $(m-k_i)$-linear maps respectively, with $k_i<m$ (i.e, multilinear maps with less than $m$ variables).
\end{definition}
The analogous definition for polynomials is: 
\begin{definition}[Schmidt Rank]\label{AR}
For a homogeneous polynomial $P$ of degree $m>1$, the Schmidt rank of $P$ as the minimal $r$  s.t:
\[
P= \sum_{i=1}^r M_iR_i
\]
where $M_i\text{ and } R_i$ are homogeneous polynomials of positive degree $<m$. 
We denote by $r(P)$ the Schmidt rank of $P$. 

For a non-homogeneous polynomial $P$ of degree $m$, let $P'$ be the $m-$homogeneous part of $P$, and we define the Schmidt rank of $P$ to be $r(P):=r(P')$.
\end{definition}
\begin{remark} \label{BRR}
Note that for a multilinear map $Q$, if $r_{pr}(Q)=s$, then $r(P_Q)\leq s$.
In addition, if P is a polynomial of degree $m$ and $r(P)=r$ we have:
\[Q_P(h_1,\ldots, h_m)=\sum_{i=1}^r \Delta h_1,\ldots \Delta h_m M_i R_i,\]
where $M_i\text{ and }R_i$ are homogeneous of degree $k_i, m-k_i$ for some $1\leq k_i<m$ respectively. Hence:
\begin{align*}
    \Delta h_1,\ldots \Delta h_m M_i R_i=\sum_{I\subset[m], |I|=k_i}\Delta_I M_i \Delta_{I^c} R_i,
\end{align*}
where:
\[\Delta_{I}M=\Delta_{h_{i_1}}\ldots\Delta_{h_{i_k}}M,\]
where $i_1,\ldots, i_k \in I.$ 
Thus:
\[Q_P(h_1,\ldots, h_m)=\sum_{i=1}^r \Delta h_1,\ldots \Delta h_m M_i R_i=\sum_{i=1}^r\sum_{I\subset[m], |I|=k_i}\Delta_I M_i \Delta_{I^c} R_i, \]
hence:
\[r_{pr}(Q_P)\leq2^m r(P).\]
\end{remark}
The notion of $r(P)$ and its connection with $r_{pr}(Q_P)$ will be important for the proof of Theorem \ref{main}.

The significance of Theorem \ref{main}  is due to the importance of polynomial sequences in the study of finite fields-geometry, and their intricate relation to Gowers Norms:
\begin{definition}
    \label{gn}(Gowers Norms)
    For a finite abelian group $G$, we define the Gowers norms of $f:G\xrightarrow{}\mC$ by:
\[\|f\|_{U^d(G)}= |\mE_{h_1,...,h_d,g\in G}\dotDelta_{h_1},...,\dotDelta_{h_d}f(g)|^{\frac{1}{2^d}}.\]
Where $\dotDelta_h f(g)=f(g+h)\overline{f(g)}$ is the multiplicative derivative in $h$.
Hence, in the setting of $\mF_q[x]$, we define the $U_d$-Gowers norm of a function $f:\mF_{q}[x]\xrightarrow{} \mC$ by:
\[\|f\|_{U^d(G_n)}= |\mE_{h_1,...,h_d,g\in G_n}\dotDelta_{h_1},...,\dotDelta_{h_d}f(g)|^{\frac{1}{2^d}}.\]
\end{definition}
It's not difficult to prove that for every integer $n$ and a character $\chi$:
\[\max_{(P,\deg(P)\leq m-1)}|\mE_{g\in G_n}f(g)\chi(P(g))|\leq ||f||_{U^{m}(G_n)}.\]
The Inverse Theorem for the Gowers Norms(\ref{ign}) implies "the converse" in the sense that the Gowers Norms of a bonded function that does not correlate with polynomials tends to zero.
Unfortunately, in section 2 we mostly deal with $
\Lambda$ (or some related functions) which is not bounded. In order to leverage Theorem \ref{ign}, we need to "replace" the boundness condition with some other more flexible condition. 

For that matter, we will define the notion of measure:
\begin{definition}
\label{measure} (Measure)
    A measure is a sequence of functions $\nu=(\nu_n)_{n\in\mN}$  s.t for every $n$, $\nu_n:G_n\rightarrow\mR_{\geq 0}$.
\end{definition}
Specifically, we are interested in measures that satisfy two technical conditions. These measures are called $M$-pseudorandom measures (sometimes, we will refer to them as measures). These conditions can be found in \cite{HLFF}.

Further essential notions for the understanding of Theorem \ref{LEP}, are the notions of complexity and normal forms.
We start by defining the complexity of an affine linear system: 
\begin{definition}(s-complexity form)\label{s-com}
Let $\Psi=(\psi_1,\ldots, \psi_t): \mF_p[x]^d\xrightarrow{}\mF_p[x]^t$ be a system of affine linear forms. For each $1\leq i\leq t$ we say that at $i$ $\Psi$ has complexity at most $s$, if we can cover the set 
$$\{\psi_j|j\neq i\}$$
by $s+1$ classes, where $\psi_i$ does not lie in the affine-linear span of any of these classes.

We say that $\Psi$ has complexity at most $s$ if $\Psi$ has complexity at most $s$ in every $1\leq i\leq t$.
\end{definition}
Another important type of system is the $s$-normal system. These systems are closely related to the concept of a system with s-complexity, but the definitions do not coincide.
\begin{definition}(s-normal form)
Let $\Psi=(\psi_1,\ldots, \psi_t): \mF_p[x]^d\xrightarrow{}\mF_p[x]^t$ be a system of affine linear forms. For each $i$, let $S_i$ be the support of $\psi_i$. We say that $\Psi$ is $s$-normal if for every $i$ there exists some $J_i\subset S_i$ such that $|J_i|\leq s+1$ and for each $i\neq k$, $J_i$ is not a subset of $S_k$.
\end{definition}
Note that if $\Psi$ is s-normal, then $\Psi$ has complexity at most $s$. For every $i$ and $e\in J_i$, we can define:
\[A_{e,i}=\{\psi_j|\psi_j(e)=0\}\]
note, that since for every $j\neq i$, $J_i$ is not a subset of $S_j$, there exists some $e\in J_i$ s.t $\psi_j\in A_{e,i}$  hence $(A_{e,i})_{e\in J_i}$ covers $(\psi_1,\ldots, \psi_t)$, where $|J_i|\leq s+1$ thus there are at most $s+1$ classes. Furthermore, $\psi(e)\neq 0$ for every $e\in J_i$, hence $\psi_i$ does not lie in the affine-linear span of any other class, namely $\Psi$'s complexity is at most $s$.

Unfortunately, the other direction isn't true. There is, though, a way to "extend" a system of complexity at most $s$ to an $s$-normal system, as was done in Lemma \ref{extension}.

Another needed crucial definition for Theorem \ref{LEP} are the constants $\{ \beta_r \}$ that appear in Theorem \ref{LEP}:
\begin{definition}(The Von-Mangold function in finite fields)
 For every $f$ monic, define:
\begin{equation*}
\Lambda_f(g)=
\begin{cases}
\frac{|f|}{\phi(f)}& gcd(f,g)=1\\
0 & else
\end{cases}
\end{equation*}
where $\phi$ is Euler's function.   
\end{definition}
Finally, the $\{ \beta_r \}$ constants are defined as follows:
\begin{definition}($\beta_f$)\label{local_factor}
For a monic $f$ we define:
   \[\beta_f=\mE_{g\in G_{\deg(f)}^d}\prod _{i\in [t]}\Lambda_f(\psi_i(g)).\] 
\end{definition}

Note that by the Chinese remainder theorem and by the multiplicativity of Euler's function $\Phi$, we have that if :
\[\beta_a=\prod_{\text{r is irr. monic divides a}}\beta_r\]
\section{Auxiliary results}
\label{auxiliary}
In this section, we present several needed auxiliary results for the proof of theorems \ref{main} and \ref{LEP}.

The first is a version of Vaughan Lemma which is a classical result in analytical number theory that has many versions. The version below appears in \cite{BHL}.

 \begin{lemma}[Vaughan's Lemma \cite{BHL}]
    \label{VL}
    Let $\Phi$ be a bounded function from $F_q[x]$ to $\mC$, let $u,v\in \mR_{>0}$ be such that  $u,v<n$ and $v\leq n-u$.\\
    Then:
    \[
    \sum_{g\in G_n}\mu(g)\Phi(g) = -T_1+T_2
    \]
    where 
    \[
    T_1=\sum_{d, \\deg(d)\leq p^{u+v}} a_d\sum_{w\in G_{n-\\deg(d)}}\Phi(dw)
    \]
    \[
    T_2=\sum_{d, p^v\leq \\deg(d)\leq p^{n-u}}b_d\sum_{w\in G_{n-\\deg(w)}}\mu(w)\Phi(dw)
    \]
    and $$\max\{|a_d|,|b_d|\}\leq \tau (d)$$\\
    such that $\tau(d)$ is the number of monic divisors of d. 
 \end{lemma}
    $T_1 \text{ and } T_2$ are called type $\mbox{I}$ and  $\mbox{II}$ respectively.\\
    
The following lemma is a corollary of Vaughan's Lemma which was proved in \cite{BHL} . Although the lemma in \cite{BHL} assumes different conditions on $k$, the proof is essentially the same. 

 \begin{lemma}[Prop 15 in \cite{BHL}]
        \label{L1}
        Let $n\in\mN$ and $\Phi$  be a bounded function from $G_n$ to $\mC$. Let $0<c<0.5$, and let $\delta > 0$ be  s.t $$|\sum_{f \in G_n}\mu(f)\Phi(f)|> q^n \delta$$ then at least one of the two following holds:
        \begin{enumerate}
            \item There exists some $k \leq 2n^{2c}$ such that:
            $$\mE_{d \in A_{k}}|\mE_{w \in G_{n-k}}\Phi(dw)|^2\geq \frac{\delta^2}{16n^5}.$$
            \item There exists some $n^{2c}\leq k \leq n-n^{2c}$ s.t:
            \[\mE_{w,w' \in G_{n-k}}\mE_{d,d'\in A_{k}}\Phi(dw)\overline{\Phi(d'w)\Phi(dw')}\Phi(d'w')\geq \frac{\delta^4}{256n^{10}}.\]
        \end{enumerate}
        Where $A_k$ is the set of monic polynomials of degree equals to $k$.
\end{lemma}
The next needed property is an implication of the linear case which was proven in \cite{BHL} and independently in \cite{Porrit}. We use this result in the proof of our main theorem as the induction base.
\begin{corollary}\label{C1}
There exist some $0<t_1<1$ and some large enough $N_1$ , such that for every prime $p>1$, $q=p^s$, an homogeneous polynomial  $l$ of degree $1$ on $F_q[x]$, a character $\chi$ of  $\mF_q[x]$ and $n>N_1$ the following holds:
\[
|\mE_ {g\in G_n}\mu(g)\chi(l(g))| \leq q^{-n^{t_1}}.
\]
\end{corollary}
One of the results that we use in this paper claims that given a "small" set of homogeneous polynomials of bounded degree and a "large" $k$, their set of degree-$k$ common zeros is relatively large. This is a classical result, which among others, appears in the following paper:
\begin{lemma}[Prop A.3 in the Appendix, \cite{PR}] \label{Band }
    Let $V$ be an  $\mF_q$ vector space, and $Pr(V)$ be its corresponding projective space.
    Let $P_1,\ldots, P_n$ be homogeneous polynomials over $V$ of degrees $\geq 1$, and let  $$D=\sum_{i=1}^n \\deg(P_i).$$
    Let us denote
    \[\tilde{Y}=\{d\in V \mid P_i(d)=0 ,1\leq i\leq n\},\]
    and let $Y$ be the corresponding set of $Pr(V)$. Then:
    \[| Y| \geq \frac{|Pr(V)|}{2q^{D+1}}.\]
    \end{lemma}   
Another rather technical needed lemma implies that if a polynomial $P$ of rank $r$ is homogeneous and linear in some variables, then $P$ can be written as a sum of "not too many" products of the form $A\cdot B$ for  $A \text{ and }B$ homogeneous polynomials of degree smaller than $r$, where $A \cdot B$ and $P$  are linear in the same variables.
\begin{lemma}
\label{LL}
Let $P(y_1,\ldots,y_n,x_1,\ldots,x_m)$ be a homogeneous polynomial of degree $d$ and of Schmidt rank $r$, such that $P$ is $x_i$ linear for every $1\leq i \leq m$. It follows that:
\[P=\sum_{j=1}^{2^m r} Q_j R_j\]
where $Q_j$ and $R_j$ are homogeneous polynomials of degree less than d, and for each i,j either $Q_j$ or $R_j$ is of degree $1$ in $x_i$ and the other is constant in $x_i$.
\end{lemma}
   \prf\
Note that $r(P)=r$, so by the definition there exist some homogeneous polynomials $\tilde {Q_i}, \tilde{R_i}$  of degree smaller than $d$, s.t $P'=\sum_{i=1}^{r} \tilde{Q_i}\tilde{R_i}$.

For every $I \subseteq [m]$, and a polynomial $M$, we denote by $M_I$, $M$'s linear component in $(x_i)_{i\in I}$ which is constant in $(x_i)_{i\in I^C}$ (might be zero).

So we have:
\[ P=P_{[m]}=(\sum_{i=1}^{r}\tilde{Q_i}\tilde{R_i})_{[m]}=\sum_{i=1}^{r}(\tilde{Q_i}\tilde{R_i})_{[m]}=\sum_{i=1}^{r}\sum_{ I\subseteq [m]}(\tilde{Q_i})_I(\tilde{R_i})_{I^c}
\]
as desired.$\square$\\

In addition, we would like to quote a classical result regarding the Schmidt rank of the restriction of polynomials to a subspace. This result appears, for example, in \cite{CO}.
\begin{lemma} \label{sub}
Let $V$ be an $\mF$ vector space, and $P:V\xrightarrow{} \mF$ be a polynomial. Let $W\subset V$ be a subspace of  $V$ with co-dimension $k$. Denote by $P_W$ the restriction of $P$ to $W$.

Then:
\[r(P_W)\geq r(P)-k.\]
\end{lemma}

Our main ingredient is an analytical way of measuring the rank of multilinear maps.
    \begin{theorem}[Bias-rank \cite{BRc}]
    \label{BR}
    Let $V$ be an $n$-dimensional vector space over $\mF_q$. Let $Q$ be an $m$-linear map over $V^m$ and $\chi_1:\mF_q\xrightarrow{} \mC$ be an additive character defined as \[\chi_1(a)=\exp(\frac{2\pi i \operatorname{Tr}(a)}{p}).\] Then:
    \begin{enumerate}
        \item  $\mE_{x \in V^m} \chi_1(Q(x)) \geq q^{-r_{pr}(Q)}  $. (In particular, it is non-negative). 
        \item If  $\mE_{x \in V^m} \chi_1(Q(x))\geq  q^{-s}  $, then $r_{pr}(Q) < C(m)s^{C(m)}$, where $C(m)$ is a constant that depends only on $m$. (I.e. bias of the exponential sum implies low rank).\\
    \end{enumerate}
    The second condition holds as well for every polynomial and its Schmidt rank, namely:
    If  $V=\mF_q^n$ and $P:V\xrightarrow{} \mF_q$ is a polynomial of degree $m<Char(\mF_q)$, then
    $|\mE_{x \in V} \chi_1(P(x))|\geq  q^{-s} $, implies that  $r(P) < C(m)s^{C(m)}$, where $C(m)$ is a constant that depends only on $m$. (Once again, it means that the bias of the exponential sum implies low rank).
    \end{theorem}
    W.l.o.g, we can assume that $C(m)>2$ for every $m$.\\
All of the results above will assist us in proving Theorem \ref{main}. The following results will be important to us when proving \ref{LEP}.\\
We start with the Inverse Gowers Norms, which is an important tool that we use extensively in the proof of Theorem \ref{LEP}.
\begin{theorem}[The inverse theorem for the Gowers norms over function fields]\label{ign}
Let a prime $p$ and $0<s<p$ be an integer, and let $0<\delta<1$. Then, there exists $\epsilon=\epsilon_{\delta,s,p}$ s.t for every finite-dimensional vector space $V$ over $\mF_p$ and every $1$-bounded function $f$, with $\|f\|_{U^{s+1}(V)}>\delta$, there exists a polynomial $P$ of degree at most $s$ s.t:
\[|\mE_{x\in V}f(x)e_p(P(-x))|>\epsilon.\]
\end{theorem}
\begin{remark}
The main theorem in \cite{IGN} enables the removal of the characteristic condition by extending the scope from classical polynomials (defined in section \ref{notations}) to include non-classical polynomials. However, as this paper focuses exclusively on classical polynomials (Since Theorem \ref{BR} is true only with characteristic conditions), there is no need to delve into the general case here 
\end{remark}
As mentioned, this theorem is one of the main theorems which we are going to use in proving \ref{LEP}. Sadly, it is limited to the case where $f$ is a bounded function. In order to generalize the result to a wider set of functions, we will utilize The Transference Principle, which was proven by  Ho{\`a}ng L{\^e} in the context of $\mF_q[x]$ (The original  $\mZ$ case was proven by Green and Tao):
\begin{theorem}
[The Transference Principle \cite{HLFF}] \label{tp}
Let $\mF_q$ be a finite field, and $\nu=(\nu_n)_{n\in \mN}$ be a $k-$pseudorandom measure on $(\mF_{q^n}\equiv G_n)_{n\in\mN}$. Then for every $\eta>0$  and every sufficiently large $N$, every function $0\leq \phi\leq \nu$ can be decomposed as $\phi=\phi_1 +\phi_2$, where $0\leq \phi_1\leq 2+\eta$  and $\phi_2$ is uniform in the sense that $||\phi_2||_{U^{K-1}}\leq \eta$, where $K=q^k$ .  
\end{theorem}
As mentioned, Theorem \ref{tp} will assist us in extending Theorem \ref{ign} to a bigger set of functions, which will be useful for us in section \ref{proof of lep}.
Another important Lemma involving functions which are bounded by a pseudorandom measure is the following: 
\begin{theorem}[Generalized von Neumann]\label{VN}
Let $L,s,t,d\in\mN$.  Then there exists $D=D(L,s,t,d)$ such that the following holds: Let $\nu$ be a $D-$pseudorandom measure, and suppose that $f_1,\ldots, f_t:G_n\xrightarrow{}\mR$ are functions with $|f_i(x)|<\nu(x)$ for every $i$. Let $\Psi=(\psi_1,\ldots,\psi_t)$ be a system of affine-linear forms from $\mF_p[x]^d$ to $\mF_p[x]^t$ in s-normal form, with $\|\Psi\|_n<L$. Then: 
\[|\mE_{g\in (G_n)^d}\prod_{i \in [t]}f_i(\psi_i(g))|\leq \min_i \|f_i\|_{U^{s+1}(G_n)}+o(1)\]
\end{theorem}
The full proof of  Theorem \ref{VN} for $\mZ$ appears in \cite{LEP0} (The generalized Von-Neumann Lemma).

Gowers and Wolf gave a similar proof of this statement for the case of $\mF_p[x]$  in \cite{WG}. Their proof should hold for general abelian groups as well, but it only refers to the case where $f_i$ are bounded. Proving the theorem for functions controlled by a pseudorandom measure is more complicated but can be done using the same techniques that were used by Tao and Green in the proof of the theorem in \cite{LEP0} .

A similar theorem for the $\mF_p[x]$ case can also be found in \cite{HLFF}.

Note that this lemma deals with the case where $\Psi$ is s-normal, while our statement regards a system of finite complexity. To make that change, we will need the following notation and lemma:
\begin{definition}(Extension)
Let $\Psi=(\psi_1,\ldots, \psi_t): \mF_p[x]^d\xrightarrow{}\mF_p[x]^t$ be a system of affine linear forms. For some $d'\geq d$ and some affine linear system $\Psi'=(\psi_1,\ldots, \psi_t): \mF_p[x]^{d'}\xrightarrow{}\mF_p[x]^t$. We say that $\Psi'$ is an extension of $\Psi$ if:
\[\Psi(\mZ^d)=\Psi'(\mZ^{d'})\] and since we can identify $\mZ^d=\mZ^d\times\{0\}^{d'-d}$ then $\Psi$ can be seen as a restriction of $\Psi'$ to $\mZ^d$.
\end{definition}\label{extension}
In their paper \cite{LEP0}, Tao and Green proved that in the case of $\mZ$ , every system \[\Psi=(\psi_1,\ldots, \psi_t): \mZ^d\xrightarrow{}\mZ^t\] of at most $s-complexity$ has an extension \[\Psi'=(\psi_1,\ldots, \psi_t): \mZ^{d'}\xrightarrow{}\mZ^t\] who is of $s-normal$ form, where $d'=O(d)$ .If in addition, $||\Psi||_N=O(1)$, then$||\Psi'||_n=O(1)$ as well. As mentioned in \cite{WG}, their proof holds in the case of $\mF_p[x]$ as well.

Following Green's and Tao's proof, we can deduce that it is enough to prove Theorem \ref{LEP} for the case when $\Psi$ is in $s$-normal form.

To prove our theorem, we need a result similar to the one proved by Green and Tao concerning $\beta_r$. Following their paper, our statement will be as follows:
\begin{lemma}[Local factor bound]\label{beta p}
Let $\Psi=(\psi_1,\ldots, \psi_t)$ be $\mF_p[x]^d$ affine linear forms. 
Assume that for each $i$, $\psi_i$ is of the form 
$$\psi_i((g_1,\ldots, g_d))=\sum_{j=1}^d L_{i,j}g_j+b_i$$
for some polynomials $L_{i,j} \text{ and } b_i\in\mF_p[x]$, and that $\psi_i$ is not constant for every $i$. Then for every irreducible monic $r$ , $$\beta_r=1+O_{d, L, t}(p^{-\deg(r)}),$$ where $$L=\max\{\deg(L_{i,j}), \deg(b_i)| 1\leq i\leq t, 1\leq j\leq d\}.$$
In addition, if  $\{\psi_i\}_{i\leq t}$ are affinely independent in pairs, then $$\beta_r=1+O_{d,L,t}(p^{-2\deg(r)}).$$        
\end{lemma}

\prf: 
First, we can assume $\deg(r)\gg L,t,d$ (meaning that $deg(r)$ is much greater than $L,t,d$) . Note that \[\Lambda_r(\psi_i(g))=0 \Leftrightarrow (\psi_i(g),r)\neq 1\]  which means that $\sum_{j=1}^d L_{i,j}g_j\equiv -b_i \text{ mod r}$. This equation has exactly $p^{\deg(r)(d-1)}$ solutions (we may assume $L_{i,d}\neq 0$, hence for every $(g_1,\ldots, g_{d-1})$ there is exactly one solution to the equation $L_{i,d}g_d\equiv -b_i+\sum_{j=1}^{d-1} L_{i,j}g_j$).
It follows that the product $\prod _{i\in [t]}\Lambda_r(\psi_i(g))$ equals to $0$ with probability $O_t(p^{-\deg(r)})$.

If the product does not equal $0$, then it equals to $$\frac{p^{\deg(r)}}{p^{\deg(r)}-1}=1+O_{t,L,d}(p^{-\deg(r)}),$$since for an irreducible polynomial $r$ the following equation holds:
\[|r|-1=p^{\deg(r)}-1=\Phi(r).\] 
Hence:
\begin{align*}
    &\beta_r=\mE_{g\in G_{\deg(r)}^d}\prod _{i\in[t]}\Lambda_r(\psi_i(g))\\
    &=(1+O_{t,L,d}(p^{-\deg(r)}))(1-O_t(p^{-\deg(r)})=1+O_{t,d,L}(p^{-\deg(r)}),
\end{align*}
as required.

Now, if for every $i\neq j$ , $\psi_i,\psi_j$ are affinely independent, then $\psi_i,\psi_j$ can be simultaneously divided by r with probability $O_{t,Ld,}(p^{-2\deg(r)})$ , and the desire bound follows by the inclusion-exclusion formula.$\square$

\section{M\"obius is orthogonal to polynomials}\label{proof of theorem 1}

In this section we prove Theorem \ref{main}. We start by generalizing a statement from \cite{BHL}.
\begin{lemma}
\label{L2}
Let $P$ be a polynomial of degree $m$ of  $\mF_{q}[x]$, with $Char(\mF_q)>m$. Let $0<c<0.5$, and let $\delta > 0$ be  such that $$|\sum_{f \in G_n}\mu(f)\chi(P(f))|> q^n \delta,$$ then at least one of the two following holds:
\begin{enumerate}
        \item There exists some $k \leq 2n^{2c}$ and a monic polynomial $d$ of degree $k$, such that the polynomial $P_d:G_{n-k} \to \mF_q$ defined by $P_d(w)=P(dw)$ has Schmidt rank at most $C(m) (\\log_q(\frac{16n^5}{\delta^2}))^{C(m)}$.
        \item There exists some $k \in [n^{2c}, n-n^{2c}] $ and at least $(\frac{\delta^4}{512n^{10}})q^{2k}$ pairs of monic polynomials $d,d'$ of degree $k$, s.t the polynomial $P_{d.d'}:G_{n-k} \to \mF_q$ defined by
        $P_{d,d'}(w)=P(dw)-P(d'w)$ has Schmidt rank at most $$C(m)(\\log_q(\frac{256n^{10}}{\delta^4}))^{C(m)}.$$
\end{enumerate}
Where $C(m)$ is the constant from Theorem \ref{BR}
Note that this lemma is an implication of lemma 14 in \cite{BHL}.
\end{lemma}
    \prf \   Denote $\Phi=\chi(P)$. By lemma \ref{L1}, at least one of the following holds:
\begin{enumerate}
    \item There exists some $k \leq 2n^{2c}$ such that:
    $$\mE_{d \in A_{k}}|\mE_{w \in G_{n-k}}\Phi(dw)|^2\geq \frac{\delta^2}{16n^5}$$
    \item There exists some $n^{2c}\leq k \leq n-n^{2c}$ s.t:
    \[\mE_{w,w' \in G_{n-k}}\mE_{d,d'\in A_{k}}\Phi(dw)\overline{\Phi(d'w)\Phi(dw')}\Phi(d'w')\geq \frac{\delta^4}{256n^{10}}.\]
\end{enumerate}
If the first condition holds, let us  denote $\delta '=\frac{\delta ^2}{16n^5}$. Then, we must have at least one polynomial $d$ of degree $k$ s.t
\[ |\mE_{w \in G_{n-k}}\Phi(dw)|\geq|\mE_{w \in G_{n-k}}\Phi(dw)|^2\geq \delta '.\]
Define 
\[P_d(w)=P(dw).\] 
By Theorem \ref{BR}, we get that
\begin{align*}
    &r(P_d)\leq C(m)(log_q(\delta '^{-1})^{C(m)}
    =C(m)log_q(\frac{16n^5}{\delta^2})^{C(m)}   
\end{align*} hence the first case was dealt with. 
Otherwise, the second condition holds. Denote $\delta '=\frac{\delta^4}{256n^{10}}$.
Then, by the triangle inequality we get:
\begin{align*}
    &\delta ' \leq \mE_{w,w' \in G_{n-k}}\mE_{d,d'\in A_{k}}\Phi(dw)\overline{\Phi(d'w)\Phi(dw')}\Phi(d'w')\\
    &=\mE_{d,d'\in A_{k}}\mE_{w,w' \in G_{n-k}}\Phi(dw)\overline{\Phi(d'w)\Phi(dw')}\Phi(d'w') \\
    &\leq \mE_{d,d'\in A_{k}}|\mE_{w,w' \in G_{n-k}}\Phi(dw)\overline{\Phi(d'w)\Phi(dw')}\Phi(d'w')|\\
    &=\mE_{d,d'\in A_{k}}|\mE_{w,w' \in G_{n-k}}\chi(P(dw)-P(d'w)\overline{\chi(P(dw')-P(d'w')}|\\
    &=\mE_{d,d'\in A_{k}}|\mE_{w \in G_{n-k}}\chi(P(dw)-P(d'w)|^2\\
    &\leq \mE_{d,d'\in A_{k}}|\mE_{w \in G_{n-k}}\chi(P(dw)-P(d'w)|
\end{align*}
Hence, there must be at least $\frac{\delta '}{2} q^{2k}$  pairs of polynomials $(d,d')\in A_k$ , s.t:
\[| \mE_{w \in G_{n-k}} \chi (P(dw)-P(d'w))|\geq \frac{\delta '}{2}.\]
For each d,d' as above, define 
\[P_{d,d'}(w)=P(dw)-P(d'w).\]
Again, by Theorem \ref{BR} we have that for each $(d,d')$ as above, 
\begin{align*}
    r(P_{d,d'})\leq C(m)log_q(\delta '^{-1})^{C(m)}
    =C(m)log_q(\frac{256n^{10}}{\delta^4})^{C(m)} ,   
\end{align*}
as desired. $\square$\\

We can now prove our main statement.
\begin{theorem*}
Let $m\in\mN$. There exist some $0<t_m<1$ and $N_m\in \mN$, such that the following holds: for every prime number $p>2m$ , $q=p^s$ ,  $P:\mF_q[x] \to \mF_q$ a polynomial form of degree $m$, any character $\chi$ of  $\mF_q$ and $n>N_m$ we have:
\[
|\mE_ {g\in G_n}\mu(g)\chi(P(g))| \leq q^{-n^{t_m}}.
\]
\end{theorem*}
\prf \ We prove the theorem by induction.
For the induction base we will prove Corollary \ref{C1}:

Note that for $t=0.5$, large enough $N$ ,any $n>N$ and $q>(\frac{3\sqrt{3}}{2})^8$:
\[4q^{\frac{3n+1}{4}}(\frac{3\sqrt{3}}{2})^n<q^{n-n^{t}}.\] 
Thus for every $q>(\frac{3\sqrt{3}}{2})^8$, $n>N$ and $\alpha\in\mT$, by the main theorem in \cite{Porrit} we get:
\[|\mE_{f\in A_n}\mu(f)e_q(\alpha f)|\leq  q^{-n^{t}}.\]

Combining this result with the linear case in \cite{BHL}  and recalling that every linear form of $G_n$ is given by $f\xrightarrow{}(\alpha f)_{-1}$ for some $\alpha\in\mT$ , we obtain the following result:

There exist some $0<t_1'<1$ and $N_1'\in \mN$, such that the following holds:
for any prime  $p$,  $q=p^s$,$l:\mF_q[x] \to \mF_q$ a polynomial form of degree $1$ , a character $\chi$ of  $\mF_q$ and $n>N_1'$ we have:
\[|\mE_ {g\in A_n}\mu(g)\chi(l(g))| \leq q^{-n^{t_1'}}.\]

Note that
\[\sum_ {g\in G_n}\mu(g)\chi(l(g))=\sum_{c \in \mF_q*}\sum_{k=0}^{n-1}\sum_{g\in A_k}\mu(g)\chi(l(cg)).\]

It follows that  we can conclude that for a small enough $0<t_1<t_1'$ and a large enough $N_1>N_1'$ the following holds:\\
for any $q=p^s$ for a prime $p$ , any  $l:\mF_q[x] \to \mF_q$ polynomial form of degree $1$, any character $\chi$ of  $\mF_q$ and $n>N_1'$ we have:
\[|\mE_ {g\in G_n}\mu(g)\chi(l(g))| \leq q^{-n^{t_1}}.\]

For the induction step, assume that there exist some $t_{m-1}, N_{m-1}$ s.t  for every polynomial of degree at most $m-1$ and every character of $\mF_q$ with $Char(\mF_q)>2(m-1)$ the statement holds .\\
Choose a small enough $0<t_m$ s.t:
\begin{enumerate}
    \item $2t_mC(m)C(2m)<t_{m-1}$ (In particular $2t_mC(m)C(2m)<0.5$.)
    \item $2t_m<t_{m-1}$ and $C(m)^3t_m<t_{m-1}$
    \item $2C(m)^3C(2m)t_m<t_{m-1}$,
\end{enumerate}
where $C(m)$ is the constant from Theorem \ref{BR}. Recall that w.l.o.g  $C(m)>2$.
Let $N_m>N_{m-1}$ be a large enough integer s.t the following inequalities hold for every $n>N_m$ and $2\leq q\in\mN$: 
\begin{enumerate}
    \item $16n^5 q^{2n^{t_m}}\leq q^{3n^{t_m}}$ 
    \item $q^{ 2(2n^{2 t_mC(m)C(2m)} +    C(m)^{C(m)+1}2^{mC(m)}3^{C(m)^2}n^{t_m C(m)^2})-n^{t_{m-1}}}\leq q^{-n^{t_m}}$
    \item $256n^{10}q^{4n^{t_m}}\leq q^{5n^{t_m}}$
    \item 
        \[4C(m)(2^{2m}r_m)^{C(m)}+n^{t_m} \leq n^{t_{m-1}}\]
    Where 
    \[r_m=C(2m) \cdot (2^{m+1}C(m)5^{C(m)}n^{t_mC(m)})^{C(2m)}\]
    \item $\frac{q^{2n^{2t_m}-n^{t_m}}}{512 n^{10}}\geq q^{n^{t_m}}$
    \item $q^{n^{t_m}-2}>512n^{10}$
    \item \[ \frac{q^{\lfloor{n^{2t_mC(m)C(2m)}}\rfloor-1}-1 }{2q^{m 2^m 2^{2m}r_m+1}}>2.\]
\end{enumerate}
For the sake of contradiction, assume that we have some $q=p^s$ with $p>2m$, $P$ a polynomial of degree $m$, $\chi$ a character of $\mF_q$ and $n>N_m$ s.t:
    \[|\mE_{g \in G_n}\mu(g)\chi(P(g))|> q^{-n^{t_m}}.\]
From Now on, we identify $P$ with it's restriction to $G_n$, i.e. by $P$  we actually mean $P|_{G_n}$.\\
Applying Lemma \ref{L2} ($\delta=q^{-n^{t_m}}$), if the first condition holds, we have some $k\leq 2n^{2 t_m c(m)C(2m)}$ and $d\in A_k$, s.t the map $P_d:w\xrightarrow[]{} P(dw)$ is of rank lower then:
\[C(m)\\log_q(16n^5q^{2n^{t_m}})^{C(m)}\leq C(m)3^{C(m)}n^{t_m C(m)}.\]
By Remark \ref{BRR}, 
\[r_{pr}(Q_{P_d})\leq 2^m r(P_d),
\]
where $Q_{P_d}(g_1,\ldots, g_m)=\Delta g_1 \ldots \Delta g_m P_d(0)$. Note that since $P_d$ is not necessarily homogeneous (P is not necessarily homogeneous), $P_{Q_{P_d}}=m!(Pd)'$, where $(Pd)'$ is the homogeneous part of $P_d$.\\
Since 
\[\Delta_{g}P_d(0)=P_d(g)-P_d(0)=P(dg)-P(0)=\Delta_{dg}P(0)\]
we conclude by Theorem \ref{BR}:
\[q^{-r_{pr}(Q_{P_d})}\leq\mE_{g\in (G_{n-k})^m} \chi_1(Q_{P_d}(g))=\mE_{g\in (dG_{n-k})^m} \chi_1(Q_{P}(g)).
\]
Since $Q_{P'}=Q_P$ where $P'$ is the homogeneous part of $P$, by Theorem \ref{BR} we can conclude that:
\begin{multline*}
    \begin{aligned}
    &r(P'|_{dG{n-k}}) \leq r_{pr}(Q_{P'}|_{dG_{n-k}})=r_{pr}(Q_P|_{dG_{n-k}})\leq C(m)r_{pr}(Q_{P_d})^{C(m)}\\&\leq
    C(m)^{C(m)+1}2^{mC(m)}3^{C(m)^2}n^{t_m C(m)^2} 
    \end{aligned}
\end{multline*}

The Schmidt rank of $P'|_{dG_{n-k}}$ is the Schmidt rank of $P'$ restricted to the subspace $dG_{n-k}$, which implies, by Lemma \ref{sub}, that \[r(P')-k\leq r(P'|_{dG{n-k}}).\] Since $k\leq 2n^{2 t_mC(m)C(2m)}$, we have:
\[r(P)=r(P')\leq  2n^{2 t_mC(m)C(2m)} +    C(m)^{C(m)+1}2^{mC(m)}3^{C(m)^2}n^{t_m C(m)^2} .  
\label{r1}\]
Now, based on our induction assumption, we state the following lemma:
\begin{lemma}
\label{IL}
Let $N_{m-1}$, $t_{m-1}$  be as in the induction assumption. Then for every $r\in \mN$, $F:\mF_q^r \xrightarrow{} \mC$ $1$-bounded function, $P_{1},\ldots ,P_{r}$ polynomials of degree $m-1$, and $n >N_{m-1}$ we have:
\[|\mE_{g\in G_n}\mu(g)F(P_{1}(g),\ldots,P_{r}(g))|\leq q^{r-n^{t_{m-1}}}.\]
\end{lemma}
This lemma is a generalization of a lemma in \cite{BHL}. Although the original lemma solely deals with the case where $P_i$ are linear, the proof is essentially the same. The full proof of this lemma is given at the end of the section.

Now, note that if a polynomial $P$ is of rank $r$ and a degree $m$, then for every $g\in\mF_q[x]$
\[P(g)=\tilde{F}(R_1(g), \ldots R_r(g), M_1(g), \ldots M_r(g))\]
where 
\[\tilde{F}(x_1,\ldots,x_r,y_1,\ldots,y_r)\vcentcolon=\sum_{i=1}^r x_i y_i\]
for some $R_i, M_i$ polynomials of degree at most $m-1$. Meaning, for every $g$:
\[\chi(P)(g)=\chi(\tilde{F}(R_1(g),\ldots, R_r(g), M_1(g),\ldots, M_r(g))).\]
Define:
\[F\vcentcolon= \chi(\tilde{F}).\]
By Lemma \ref{IL}, we have:
\begin{align*}
 &|\mE_{g\in G_n}\mu(g)\chi(P(g))|\\
 &=|\mE_{g\in G_n}\mu(g)F(R_1(g),\ldots,R_r(g), M_1(g)\ldots, M_r(g))|\leq q^{2r-n^{t_{m-1}}}  
\end{align*}
By \ref{r1}, we conclude:
\begin{align*}
 &|\mE_{g\in G_n}\mu(g)\chi(P(g))|\\
 &\leq q^{ 2(2n^{2 t_mC(m)C(2m)} +    C(m)^{C(m)+1}2^{mC(m)}3^{C(m)^2}n^{t_m C(m)^2})-n^{t_{m-1}}}\leq q^{-n^{t_m}},   
\end{align*}
which contradicts our first assumption.

Otherwise, the second condition holds, meaning that there exist some \\$k \in [n^{2t_mC(m)C(2m)}, n-n^{2t_m C(m)C(2m)}] $ and at least $\frac{q^{2k-4n^{t_m}}}{512n^{10}}$ \\pairs $(d,d')\in A_k^2$ s.t: 
\[r(P_{d,d'})\leq C(m)\\log_q(256n^{10}q^{4n^{t_m}})^{C(m)}\leq C(m)5^{C(m)}n^{t_m C(m)},\]
where 
\[P_{d,d}(w)=P(dw)-P(d'w).\]
Denote
\[Y=\{ (d,d')\in G_{k+1}^2| r(P_{d,d'})\leq C(m)5^{C(m)}n^{t_m C(m)}\}.\]
Since $Y$ is of cardinality $\geq \frac{q^{2k-4n^{t_m}}}{512n^{10}}>1$, Y is not empty. 
\\For each (d,d')$\in G_{k+1}^2$, denote the $m$-linear map $Q_{d,d'}$ as:
\[Q_{d,d'}(g_{1},\ldots,g_{m})=Q(dg_{1},\ldots,dg_{m})-Q(d'g_{1},\ldots,d'g_{m})\]
Where Q=$\frac{Q_{P'}}{m!}$, i.e \[Q(g,\ldots,g)=P'(g)\] for every g, and \[Q_{d,d'}(g,\ldots,g)=P'(dg)-P'(d'g)=P'_{d,d'}(g).\] 
Where $P'$ is the homogeneous part of $P$. \\
Note that we use the fact that $m<Char(\mF_q)$, since we can not properly define $Q$ if $m\geq Char(\mF_q)$.\\
By Remark \ref{BRR}, for each $(d,d')\in Y$: 
\begin{align*}
    &r_{pr}(Q_{d,d'})\leq 2^m r(P'_{d,d'})\leq 
    2^m r((P_{d,d'})')=2^m r(P_{d,d'})\leq
    2^mC(m)5^{C(m)}n^{t_mC(m)}
\end{align*}
And since for each $(d,d')\in G_{k+1}^2$ the map $Q_{d,d'}$ is multilinear, we have (by Theorem \ref{BR}):
\begin{align*}
    &\mE_{(x_1,\ldots,x_m)\in G_{n-k}^m}\chi_1(Q_{d,d'}(x_1,\ldots,x_m))\geq 0
\end{align*}
So we have:
\begin{align*}
     &\mE_{(d,d', g_1,\ldots,g_m)\in G_{k+1}^2 \times G_{n-k}^m}\chi_1(Q_{d,d'}(g_1, \ldots, g_m))\\
    &= \mE_{(d,d')\in G_{k+1}^2}\mE_{(g_1,\ldots,g_m)\in G_{n-k}^m}\chi_1(Q_{d,d'}(g_1, \ldots, g_m))\\
    &=\frac{1}{q^{2k+2}}\sum_{(d,d')\in G_{k+1}^2}\mE_{(g_1,\ldots,g_m)\in G_{n-k}^m}\chi_1(Q_{d,d'}(g_1, \ldots, g_m))\\
    &\geq\frac{1}{q^{2k+2}}\sum_{(d,d')\in Y}\mE_{(g_1,\ldots,g_m)\in G_{n-k}^m}\chi_1(Q_{d,d'}(g_1, \ldots, g_m))\\
    &\geq\frac{1}{q^{2k+2}}\frac{q^{2k-4n^{t_m}}}{512n^{10}} q^{-r_{pr}(Q_{d,d'})}
    \geq q^{-2r_{pr}(Q_{d,d'})} \\
    &\geq q^{-2^{m+1}C(m)5^{C(m)}n^{t_mC(m)}}   
\end{align*}

(w.l.o.g $C(m)>2$).\\
Denote  $R:G_{k+1}^2\times G_{n-k}^m\xrightarrow{} \mF_q$ by:
\[R(d,d',g_{1},\ldots,g_{m})=Q_{d,d'}(g_1,\ldots,g_m)=Q(dg_{1},\ldots,dg_{m})-Q(d'g_{1},\ldots,d'g_{m}).\]
Note that $R$ is a homogeneous polynomial of degree $2m$. Since $p>2m$, by Theorem \ref{BR}
\[r(R)\leq C(2m) \cdot (2^{m+1}C(m)5^{C(m)}n^{t_mC(m)})^{C(2m)}=r_m\]
Define $R'(d,g)=P(dg)$ and $Q'$ to be the $2m$-linear map corresponding to $R'$. \\

Since $R'$ is the  restriction of $R$ to the subspace $G_k\times 0\times V$, where $V=\{(x_1,\ldots, x_m)\in G_{n-k}^m| g_1=\ldots=g_m\},$ we have:
\[r(Q')\leq 2^{2m} r(R')\leq 2^{2m}r(R)=2^{2m}r_m.\]
Note that 
\[Q'(d_1,\ldots, d_m, g_1,\ldots, g_m)=Q(d_1g_1,\ldots, d_m g_m).\]
Note that either $k\leq n-k$ or $k> n-k$. Assume that $k\leq n-k$.\\
Observe:
\[Q'=\sum_{i=1}^{r(Q')} M_i T_i\] where $M_i, T_i$ are $l_i-, (2m-l_i)-$linear maps for every $i$.\\
Namely, each term is either of type:
\begin{align}\label{1}
    &S(d_{i_1},\ldots, d_{i_m})H(d_{j_1},\ldots, d_{j_m},g_{1},\ldots,g_{m}) 
\end{align}
where $\{i_1,\ldots, i_m\} \cup \{j_1,\ldots, j_m\}=[m]$.\\
Or of type: $T\cdot A$ where $T$ and $A$ are multilinear maps that have less than $m$ variables of $G_{n-k}.$ Hence $T\cdot A$ is linear in $g_1,..., g_m$ and for each $(d_1,\ldots, d_m)$, the map $  (T\cdot A)(d_1,\ldots, d_m,\cdot)$ is $m$-linear map of rank 1.\\
Define $J=\{0, \lfloor k-n^{2t_mC(m)C(2m)}\rfloor \}$.\\ 
For each $(j_1,\ldots, j_m)\in J^m$, and $S_i, H_i$ of type \ref{1}, we can define a polynomial \[S_i^{j_1,\ldots, j_m}(d)=S_i(x^{j_1}d,\ldots,x^{j_m}d)\]of degree $m$ from $G_{\lfloor n^{2t_mC(m)C(2m)}\rfloor}$ to $\mF_q$.\\
Let \[A=\{d\in G_{\lfloor{n^{2t_mC(m)C(2m)}}\rfloor}| S_{i}^{j_1,\ldots, j_m}(d)=0,  \text{ for every i and } \{j_1,\ldots, j_m\}\in J^m\}.\]

Note that since $\deg(S_i^{j_1,\ldots, j_m})\leq m$, and since we have at most $r_{pr}(Q')$ polynomials $S_i-s$ and  $|J^m| = 2^m$, we have: 
\[D=\sum_{i,\{j_1,\ldots, j_m\}} \deg(S_{i}^{j_1,\ldots, j_m})\leq m\cdot 2^m r_{pr}(Q').\]

By Theorem \ref{PR}, \[|Pr(A)|\geq\frac{|Pr(G_{\lfloor{n^{2t_mC(m)C(2m)}}\rfloor})|}{2q^{D+1}}\geq \frac{{q^{\lfloor n^{2t_mC(m)C(2m)}\rfloor-1}-1} }{2q^{m 2^m 2^{2m} r(R)+1}}.\]

Hence, it has a common zero $w\in G_{\lfloor n^{2t_m C(m)C(2m)} \rfloor} \backslash \{0\}$.\\
For each $(j_1,\ldots, j_m)\in J^m$ we can define: 
\[Q'_{(j_1,\ldots, j_m)}(g_1,\ldots, g_m)=Q'(x^{j_1}w,\ldots, x^{j_m}w, g_1,\ldots,g_m)=Q(x^{j_1}wg_1,\ldots, x^{j_m}wg_m).\]
Since $w$ is a common zero of all $S_i^{j_1,\ldots, j_m}$, we get that for every $(j_1,\ldots, j_m)\in J^m$  $Q'_{(j_1,\ldots, j_m)}$ is an $m-linear$ map of $G_{n-k}^m$ and its rank is at most $r(Q')$ .\\
We now observe $Q$ in the subspace $wG_{\lfloor n-n^{2t_mC(m)C(2m)}\rfloor}^m$.
Note that each $g$ in this space can be uniquely written as a sum of polynomials \[v\in wG_{\lfloor k-n^{2t_mC(m)C(2m)}\rfloor}\] and \[u\in w\cdot sapn\{x^{\lfloor k-n^{2t_mC(m)C(2m)}\rfloor},\ldots, x^{\lfloor n-n^{2t_mC(m)C(2m)}\rfloor} \}=wx^{\lfloor k-n^{2t_mC(m)C(2m)}\rfloor}G_{n- k},\] Hence:
\begin{align*}\label{calc} (2)
    &\sum_{(g_1, \ldots, g_m)\in wG_{\lfloor n-n^{2t_mC(m)C(2m)}\rfloor}^m} e_q(Q(g_1,\ldots, g_m))=\\
    &\sum_{(v_1,\ldots, v_m)\in wG_{\lfloor k-n^{2t_mC(m)C(2m)}\rfloor}^m}\sum_{(u_1,\ldots, u_m)\in wx^{\lfloor k-n^{2t_mC(m)C(2m)}\rfloor}G_{n- k}^m}e_q(Q(v_1+u_1,\ldots, v_m+u_m))=\\
    &\sum_{(v_1,\ldots, v_m)\in wG_{\lfloor k-n^{2t_mC(m)C(2m)}\rfloor}^m}\sum_{(u_1,\ldots, u_m)\in wx^{\lfloor k-n^{2t_mC(m)C(2m)}\rfloor}G_{n- k}^m}\sum_{\substack{(g_1,\ldots, g_m),\\ g_i\in\{u_i, v_i\}}}e_q(Q(g_1,\ldots, g_m))=\\
    &\sum_{\substack{(v_1,\ldots, v_m)\in\\G_{\lfloor k-n^{2t_mC(m)C(2m)}\rfloor}^m}}\sum_{(u_1,\ldots, u_m)\in G_{n- k}^m}\sum_{, (j_1,\ldots,j_m)\in J^m}\sum_{\substack{(g_1,\ldots, g_m), g_i\in\{u_i, vi\} \\g_i=\begin{cases}
        v_i& j_i=0\\
        u_i & else 
    \end{cases}}}(e_q(Q(wx^{j_1}g_1,\ldots, wx^{j_m}g_m))=\\
    &\sum_{\substack{(v_1,\ldots, v_m)\in\\G_{\lfloor k-n^{2t_mC(m)C(2m)}\rfloor}^m}}\sum_{(u_1,\ldots, u_m)\in G_{n- k}^m}\sum_{, (j_1,\ldots,j_m)\in J^m}\sum_{\substack{(g_1,\ldots, g_m), g_i\in\{u_i, vi\} \\g_i=\begin{cases}
        v_i& j_i=0\\
        u_i & else 
    \end{cases}}}e_q(Q'_{j_1,\ldots, j_m}(g_1,\ldots, g_m))
\end{align*}
Now, for each $(j_1,\ldots, j_m)$ , $Q'_{j_1,\ldots,j_m}$ depending only on $V_{j_1,\ldots, j_m}=\prod_{i=1}^m V_i$, where:
\[V_i=
\begin{cases}
G_{\lfloor k-n^{2t_mC(m)C(2m)}\rfloor} & j_i=0\\
G_{n-k} & \text{else}
\end{cases}.\]
Assume that $j_i=0$ s-times. Then $|V_{j_1,\ldots, j_m}|=q^{s\lfloor k-n^{2t_mC(m)C(2m)}\rfloor}q^{(m-s)(n-k)}$. \\
Hence, for each $(g_1,\ldots, g_m)\in V_{j_1,\ldots, j_m}$, the expression $Q'_{(j_1,\ldots, j_m)}(g_1,\ldots, g_m)$ appears in (2) $$\frac{q^{m(\lfloor n-n^{2t_mC(m)C(2m)}\rfloor)}}{q^{s\lfloor k-n^{2t_mC(m)C(2m)}\rfloor}q^{(m-s)(n-k)}}\text{-times.}$$\\
Also, we may notice that for each $(j_1,\ldots, j_m)$ $V_{j_1,\ldots, j_m}$ is a subspace of $G_{n-k}^m$ (because $n-k\geq k$),  hence we can bound the rank of $Q'_{j_1,\ldots, j_m}$ by $r(Q')$.\\
Hence:
\begin{align*}
    &\sum_{\Bar{g}\in (wG_{\lfloor n-n^{2t_mC(m)C(2m)}\rfloor}^m)}   e_q(Q_P(\Bar{g}))=\\
    &\sum_{s=0}^m \sum_{\substack{\Bar{j}\in J^m\text{ s.t }\\|\{i\in \Bar{j}|i=0\}|=s}}\frac{q^{m(\lfloor n-n^{2t_mC(m)C(2m)}\rfloor)}}{q^{s\lfloor k-n^{2t_mC(m)C(2m)}\rfloor}q^{(m-s)(n-k)}}\sum_{\Bar{g}\in V_{\Bar{j}}}e_q(Q'_{\Bar{j}}(\Bar{g}))
\end{align*}

By dividing by ${q^{m(\lfloor n-n^{2t_mC(m)C(2m)}\rfloor)}}$, we get:
\begin{align*}
    &\mE_{(g_1, \ldots, g_m)\in (wG_{\lfloor n-n^{2t_mC(m)C(2m)}\rfloor}^m)} e_q(Q_P(g_1,\ldots, g_m))=\\
    & \sum_{s=0}^m \sum_{\substack{(j_1,\ldots, j_m)\in J^m\text{ s.t } \\ |\{i|j_i=0\}|=s }}\mE_{(g_1,\ldots, g_m)\in V_{j_1,\ldots, j_m}}e_q(Q'_{j_1,\ldots, j_m}(g_1,\ldots, g_m))\geq\\
    &\sum_{\substack{(j_1,\ldots, j_m)\in J^m}}q^{-r(Q'_{j_1,\ldots, j_m})}\geq 2^mq^{-r(Q')}\geq q^{-r(Q')}
\end{align*}
By using Bias-Rank, we get that:
\[r(Q|_{(wG_{\lfloor n-n^{2t_mC(m)C(2m)}\rfloor})^m})\leq C(m)r(Q')^{C(m)}\]
and Hence:
\begin{align*}
    &r(P')\leq r(P'|_{wG_{\lfloor n-n^{2t_m C(m)C(2m)}\rfloor}})+n^{2t_mC(m)C(2m)}\\
    &\leq r(Q_P|_{wG_{\lfloor n- n^{t_mC(m)C(2m)}\rfloor}^m})+n^{2t_mC(m)C(2m)}\\
    &\leq C(m)r(Q')^{C(m)}+n^{2t_mC(m)C(2m)}\leq 2C(m)r(Q')^{C(m)}.
\end{align*}
Note that the case where $k>n-k$  is symmetric with respect to $k, n-k$, Hence we can repeat the calculation by replacing $k$ with $n-k$ and $(g_1,\ldots, g_m)$ with $(d_1,\ldots, d_m)$, and we obtain the same bound on $r(P')$.
\\
By the same arguments of the first case, we are done. $\square$\\
The only remaining thing is the proof of  Lemma \ref{IL}.
As mentioned before, the lemma was already proven in \cite{BHL} in the case of a polynomial of degree $1$, and the proof of the general case is essentially the same.

{\em Proof of Lemma \ref{IL}}. By the induction assumption,
for every polynomial of degree $m-1$, $P$, and $n>N_{m-1}$ the following holds: \[
|\mE_{g\in G_{n}}\mu(g)\chi(P(g))|\leq q^{-n^{t_{m-1}}}.
\]

For every $a=(a_1,\ldots,a_r)\in \mF_{q}^r$, denote 
\[
V_a=\{g\in G_n| P_{i}(g)=a_i \ \text{for every}\  i\}.
\]
Then:
\[
\sum_{g\in G_n}\mu(g)F(P_{1}(g),\ldots,P_{r}(g))=\sum_{a\in F_{p}^r}F(a_1,\ldots,a_r)\sum_{g\in V_a}\mu(g).
\]
Note that 
\[
1_{V_a}(g)=\mE_{(\chi_1,\ldots,\chi_r)=\chi\in \hat{F_{q}^r}}\prod_{i\in[r]} \chi_i(P_i(g)-a_i),
\]
hence 
\[
\sum_{g\in V_{a}}\mu(g)=\mE_{(\chi_1,\ldots,\chi_r)=\chi \in \hat{F_{p}^r}}\prod_{i \in [r]}\chi_i(-a_i)\sum_{g \in G_n}\mu(g)\prod_{i\in [r]}\chi_i(P_i(g)).
\]
By the triangle inequality, 
\[
|\sum_{g\in V_{a}}\mu(g)|\leq \max_{\chi}|\sum_{g\in G_n}\mu(g)\prod_{i \in [r]}\chi_i(P_i(g))|.
\] 
Note that each $\chi_i(x)=e_q (s_i x)$ for some $s_i$, so we have :
\[
\prod_{i\in [r]}\chi_i(P_i(g))=e_q(\sum_{i\in[r]} s_i P_i(g))
\] 
when $\sum_{i\in[r]} s_i P_i(g)$ is a polynomial of degree $m-1$.
So we have that for any $n >N_{m-1}$:
\[
|\sum_{g\in V_{a}}\mu(g)|\leq q^{n-n^{t_{m-1}}}.
\]
It follows that
\begin{align*}
    &|\mE_{g\in G_n}\mu(g)F(P_{1}(g),\ldots,P_{r}(g))|\\
    &=q^{-n}|\sum_{a\in F_{p}^r}F(a_1,\ldots,a_r)||\sum_{g \in V_a}\mu(g)|\leq p^{r-n^{t_{m-1}}}    
\end{align*}
as desired.$\square$

\section{Linear Equations in Primes}\label{proof of lep}
In this section we will prove Theorem \ref{LEP}. We will follow the proof that Green and Tao outlined in \cite{LEP0} to establish Theorem \ref{LEP}.

First, define:
\[\Lambda'(g)=\begin{cases}
    \deg(g) &\text{ if g is prime}\\
    0 &\text{ else}
\end{cases}\]
We may use $\Lambda$ and $\Lambda'$ are interchangeable since $\Lambda'$ differs from $\Lambda$ on a negligible set of prime powers. Ultimately, we prove Theorem \ref{LEP} for $\Lambda'$.

Our strategy for proving this version of the theorem involves the usage of the "W-trick" in order to replace the von Mangoldt function $\Lambda$ with a similar, more easily estimable function $\Lambda'_{W,b}$, which we will formally define later.

Specifically, we construct some pseudorandom measure $\nu$ that controls $\Lambda'_{W,b}-1$, and then use Theorem \ref{VN} combined with theorem \ref{main}.

We start with some function $w:\mN\xrightarrow{}\mR$ which tends slowly to infinity. Specifically, we can choose: \[w(n)=\frac{1}{2}\log(\log(n))\]
For some $N\in\mN$ we define \[W=W(w(N))=\prod_{r \text{ irreducible monic, } \deg(r)\leq w(N)}r.\]
By using the prime number theorem we deduce:
\begin{align*}
    &\deg(W)\leq |\{\text{ r monic irreducible, } \deg(r)\leq w(N)\}|w(N)\\
    &\leq p^{w(N)}w(N)\leq(\log (N))^{0.75}.
\end{align*}
For the purpose of this paper, it is enough to take any $W=o(\log(n)).$\\
From now on, we mainly deal with a function which is related to $\Lambda$ defined as the following:
\begin{definition}($\Lambda_{W,b}$)\label{LWb}
Take some $b$ of degree lower than $\deg(W)$, and $gcd(W,b)=1$. We define:
\[\Lambda_{W,b}(g)=\Lambda(Wg+b).\]
\end{definition}
We can now state the following theorem:
\begin{theorem}\label{bL}
For a polynomial of degree $m$ less than $\frac{p}{2}$, and $W,b$ as in definition \ref{LWb}.  The following holds:
 \[|\mE_{g\in G_n} (\frac{\Phi(W)}{|W|}\Lambda_{W,b}(g)-1)\chi(P(g))|=o_{m}(1)\]
 Where $\Phi$ is Euler's function defined over $\mF_p[x]$.
\end{theorem}
\prf:\\
Note that:
\[\Lambda(g)=-\sum_{d|g}\mu(d)\deg(d).\]
We now perform a smooth splitting of $Id$ to $Id=Id^{\sharp}+Id^{\flat}$, where $Id^{\sharp}$ is supported on $\|x\|\leq 1$, and $Id^{\flat}$ is supported on $\|x\|\geq 0.5$, while both are bounded by $Id$.
In addition, we want that $Id^{\sharp}$ will satisfy that $\int_0^{\infty}(Id^{{\sharp}'})^2=1$. We may also assume that $Id^{\sharp}$ vanishes on a small neighborhood of $0$ (In \cite{YT} one can find the full details of the splitting).
This induced a splitting of $\Lambda$ to $\Lambda=\Lambda^{\sharp}+\Lambda^{\flat}$ in the following way:
\begin{align*}
    &\Lambda^{\sharp}(g)=-R\sum_{d|g}\mu(d)Id^{\sharp}(\frac{\deg(d)}{R})\\
    &\Lambda^{\flat}(g)=-R\sum_{d|g}\mu(d)Id^{\flat}(\frac{\deg(d)}{R})
\end{align*}
for $R=sN$ (Where $=s>0$ is sufficiently small).\\
Hence it is sufficient to show that:
\begin{align}
 &|\mE_{g\in G_n}(\frac{\Phi(W)}{|W|} \Lambda^{\sharp}(Wg+b)-1)\chi(P(g))|=o_{m}(1)\label{sharp}\\
 & |\mE_{g\in G_n} \frac{\Phi(W)}{|W|}\Lambda^{\flat}(Wg+b)\chi(P(g))|=o_{m}(1)\label{flat}
\end{align}
Both cases, \ref{sharp} and \ref{flat}, were proven By Y. Cohen (\cite{YT}, Propositions 6.3 and 6.4) for a sufficiently small $s>0$. His proof is based on  Ho{\`a}ng L{\^e} paper \cite{HLFF} and on Theorem \ref{main}.$\square$

Ideally, we want to use the Inverse Gowers Norms (Theorem \ref{ign}) to prove that $\frac{\Phi(W)}{|W|}\Lambda_{W,b}-1$ has small Gowers Norms. Unfortunately,  $\frac{\Phi(W)}{|W|}\Lambda_{W,b}-1$ is not bounded, hence we cannot directly apply Theorem \ref{ign} .

To overcome this obstacle we will use Theorem \ref{tp} and deduce the following relative inverse theorem:
\begin{theorem}[Relative Inverse Gowers Norms Theorem]\label{rign}
Let $0<s<p$ be an integer, and let $0<\delta<1$. Then there exists $\epsilon=\epsilon_{\delta,s,p}$ and $D=D_s\in\mN$ s.t for every finite-dimensional vector space $V$ over $\mF_p$ and every function $f$ bounded pointwise by a D-pseudorandom measure $\nu$, with $||f||_{U^{s+1}(V)}>\delta$, there exists a polynomial $P$ of degree at most $s$ s.t:
\[|\mE_{x\in V}f(x)e_p(-P(x))|>\epsilon.\]
    
\end{theorem}
\prf:
Assume $f$ satisfies the conditions above.  By Theorem \ref{tp} , we can decompose $f$ to  $f=f_1+f_2$  , where $f_1$ is bounded, and $f_2$  has a small $q^D-1$ Gowers norm. By choosing $D$ sufficiently large, and by the monotony of the Gowers norms, we deduce $f_2$  has a small $s+1$ Gowers norm.

By the triangle inequality, we conclude there exists some $\delta'>0$ such that:
\[||f_1||_{U^{s+1}}>\delta'.\]
Since $f_1$ is bounded, we can apply Theorem \ref{ign} to $f_1$ and find some polynomial $P$ of degree at most $s$, and $\epsilon'>0$ such that the following holds:
\[|\mE_{x\in V}f_1(x)e_p(-P(x))|>\epsilon'.\]
By the monotony of the Gowers Norms, we deduce that for every polynomial $P$ of degree at most $s$:
\[|\mE_{x\in V}f_2(x)e_p(-P(x))|\leq ||f_2 \cdot e_p(-P)||_{U^{s+1}}=||f_2||_{U^{s+1}}\]
Where the last equality is true for every $f$ and every polynomial of degree at most $s$, and is cited, for instance, in \cite{IGN}.

Using the triangle inequality again, we claim that there exists some $\epsilon>0$ such that:
\[|\mE_{x\in V}f(x)e_p(-P(x))|>\epsilon\]
as required. $\square$

Fortunately, we can bound $\Lambda'_{W,b}-1$ by a pseudorandom measure.\\
First, since the difference between $\Lambda_{W,b}$ and $\Lambda'_{w,b}$ is negligible, we may assume that Theorem \ref{bL} holds also for $\Lambda'_{W,b}$.
\begin{claim}\label{dpm}(Domination by a Pseudorandom measure)
    For every $D\in\mN$, there exists an $D$-pseudorandom measure $\nu$ and a positive small $c=c_{s,t,d,L}$ s.t for every $g$ of large enough degree: $$c|( \Lambda'_{W,b}-1)(g)|<\nu(g)$$

\end{claim}
\prf:
Take some $m$.
Define:
\[\Lambda_{\chi, R}(g)=\sum_{d|g}\mu(d)\chi(\frac{\deg(d)}{R})\]
Where $\chi:\mR\rightarrow\mR$ is a smooth function supported on $[-1,1]$, such that $\chi(0)>0$ and $\int_{0}^{\infty}(\chi'(x))^2dx=1$.

First let us prove that for every  $g$ of degree at least $\beta N$ where $1>\beta>s$ (Where $R=sN$) :
\[\frac{\Phi(W)}{|W|}\Lambda'_{W,b}(g)\ll\frac{\Phi(W)}{|W|} R\Lambda_{\chi, R}(Wg+b)^2\]
The left-hand side is non-zero only if $Wg+b_i$ is irreducible. In that case, it is equal to $\frac{\Phi(W)}{|W
|}sN<\frac{\Phi(W)}{|W
|}\deg(Wg+b_i)<\frac{\Phi(W)}{|W|}(N+deg(W))$ and the right-hand side equals to $\frac{\Phi(W)}{|W|}|\chi(0)|^2R$.
Since $R=O(N)$ the claim follows.\\
Hence:
\[\frac{\Phi(W)}{|W|}\Lambda'_{W,b}(g)+1\ll\frac{1}{2}\frac{\Phi(W)}{|W|} R\Lambda_{\chi, R}(Wg+b)^2+\frac{1}{2}\label{*}(*).\]
 Ho{\`a}ng L{\^e} Proved in \cite{HLFF} that for a small enough  $s=s_D$, the right side is a $D$-pseudorandom measure .\\
By the inequality, we can find a small positive $c$ such that:
\[c|\frac{\Phi(W)}{|W|}\Lambda'_{W,b}(g)-1|<\frac{1}{2}\frac{\Phi(W)}{|W|} R\Lambda_{\chi, R}(Wg+b)^2+\frac{1}{2}\]
as desired.$\square$\\
Now, by using the Relative Inverse Gowers Norms (Theorem \ref{rign}) we can deduce:
\begin{corollary}\label{sgn}
    For any $s\leq\frac{p}{2}$: \[||\frac{\Phi(W)}{|W|}\Lambda'_{W,b}-1||_{U^{s}(G_n)}=o_s(1)\]
\end{corollary}
Hence, we can now deduce the following Theorem:
\begin{theorem}\label{LIP}
For each $(b_1,\ldots, b_t)$ that is coprime to $W$ in each coordinate, and any $s-$normal system $\Psi=(\psi_1,\ldots, \psi_t)$ of affine-linear forms for some $s<\frac{p}{2}$, we have:
\[\mE_{g\in G_n}\prod_{i
\in[t]}(\frac{\Phi(W)}{|W|}\Lambda'_{W,b_i}(\Psi_i(g))-1)=o_s(1)\]   
\end{theorem}
\prf:\\
We start by using \ref{dpm} to prove that there exists an $D$-pseudorandom measure $\nu$ that dominates $\frac{\Phi(W)}{|W|}\Lambda'_{W,b_i}$ for each $i$. By Lemma \ref{dpm} we have for each $i$:
\[\frac{\Phi(W)}{|W|}\Lambda'_{W,b_i}(g)\ll\frac{\Phi(W)}{|W|} R\Lambda_{\chi, R}(Wg+b_i)^2\]Hence:
\[1+\sum_{i\in [t] }\frac{\Phi(W)}{|W|}\Lambda'_{W,b_i}(g)\ll\frac{1}{2}+\frac{1}{2}\mE_{i\in[t]}\frac{\Phi(W)}{|W|} R\Lambda_{\chi, R}(Wg+b_i)^2=\nu.\]
Where $\nu$ is a pseudorandom measure. \\
In particular, there exists some small positive $c$  s.t for each $i$, \[c|\frac{\phi(W)}{|W|}\Lambda'_{W,b_i}-1|<\nu.\label{bbm}\]
Applying Lemma \ref{VN} , by choosing $f_i=c(\Lambda'_{W,b_i}-1)$ and using claim \ref{sgn}, we get:
\[\mE_{g\in G_n}\prod_{i
\in[t]}(\frac{\Phi(W)}{|W|}\Lambda'_{W,b_i}(\Psi_i(g))-1)=o(1)\]
as desired $\square$.\\

Now, note that any subsystem of $\Psi$ is still $s$-normal. Therefore, we can rewrite $\Lambda_{W,b_i}$ as $((\Lambda_{W,b_i}-1)+1)$. By applying Theorem \ref{LIP} repeatedly, we obtain:
\begin{theorem}
    For each $(b_1,\ldots, b_t)$ that is coprime to $W$ in each coordinate, and any $s-$normal system $\Psi=(\psi_1,\ldots, \psi_t)$ of affine-linear forms for some $s<\frac{p}{2}$, we have:
    \[\mE_{g\in G_{n}^d}(\prod_{i \in [t]}\frac{\Phi(W)}{|W|}\Lambda '_{W,b_i}(\psi_i(g))-1)=o(1)\]
\end{theorem}
We now have everything we need to prove our main theorem, Theorem \ref{LEP}. 

Proof Of Theorem \ref{LEP}:\\
\[\sum_{g\in G_n^d}\prod_{i\in[t]}\Lambda'(\psi_i(g))=\sum_{a\in G_{\deg(W)}^d}\sum_{g\in G_{n-\deg(W)}^d}\prod_{i\in [t]}\Lambda'(\psi_i(Wg+a))\]
Define:
\[A=\{a\in G_{\deg(W)}^d| gcd(W, \psi_i(a))=1 \text{ for every i}\}\]
Then for every $a\notin A$, there exists some $i$ s.t:
\[\psi_i(Wg+a)=\Dot{\psi_i}(Wg)+\Dot{\psi_i}(a)+\psi_i(0)=W\Dot{\psi_i}(g)+\psi_i(a)\]
 hence $gcd(\psi_i(Wg+a), W)\neq 1$. Hnece there exists some irreducible monic polynomial $r$ of degree$\leq w(n)$, s.t $r$ divides $\psi_i(Wg+a)$, which implies that $\psi_i(Wg+a)$ is irreducible if and only if $\psi_i(Wg+a)\in r\mF_p$. In that case, $\Lambda'(\psi_i(Wg+a))=\deg(r)\leq w(n)$.\\
Note that w.l.o.g,  for each $i$  $L_{i,d}\neq 0$ , and for  $(g_1,\ldots, g_{d-1})\in G_{n-\deg(W)}^{d-1}$ and $\alpha\in\mF_p$, $\psi_i(g_1,\ldots, g_{d-1},g)=\alpha r$ has at most one solution. Hence, the set: 
\[\{g\in G_{n-\deg(W)}^d| \text{ for each i } \psi_i(Wg+a)\in r\mF_p\}\]
is of size less or equal to: 
\[p|G_{n-\deg(W)}^{d-1}|=p(p^{(n-\deg(W))(d-1)})=p^{d(n-\deg(W))-(n-\deg(W))+1}\]
and hence:
\begin{align*}
    &\sum_{a\notin A}\sum_{g\in G_{n-\deg(W)}^d}\prod_{i\in [t]}\Lambda'(\psi_i(Wg+a))\\
    &\leq \sum_{a\notin A}\sum_{r \text{ irr. } r|W}\sum_{g\in G_{n-\deg(W)}^d,\text{ for every i } \psi_i(Wg+a)\in r\mF_p}w(n)\\
    &\leq p^{d\deg(W)} p^{w(n)} p^{d(n-\deg(W))-(n-\deg(W))+1}w(n)=o(p^{nd}).
\end{align*}

For every $i$, if $a\in A$, then we can write:
\[\psi_i(Wg+a)=W\Tilde{\psi}_{i,a}(g)+b_i(a)\]
where $b_i(a)$ is co-prime to $W$, $0<\deg(b_i(a))<\deg(W)$ and $b_i(x), c_i(a)$ satisfy 
\[\psi_i(a)=Wc_i(a)+b_i(a)\]
and  
\[\Tilde{\psi}_{i,a}(g)=\Dot{\psi_i}(g)+c_i(a).\]
Note that $\Tilde{\Psi}=(\Tilde{\psi}_{1,a},\ldots,\Tilde{\psi}_{t,a})$ has the same complexity as $\Psi$.\\
Now, note that:
\[\Lambda'(\psi(Wg+a))=\frac{|W|}{\phi(W)}\frac{\phi(W)}{|W|}\Lambda'_{W, b_i(a)}(\Tilde{\psi}_{i,a}(g))\]
meaning:
\begin{align*}
    &\sum_{g\in G_n^d}\prod_{i\in[t]}\Lambda'(\psi_i(g))=\sum_{a\in A}(\frac{|W|}{\phi(W)})^t\sum_{g\in G_{n-\deg(W)}^d}\prod_{i\in [t]}\frac{\phi(W)}{|W|}\Lambda'_{W, b_i(a)}(\Tilde{\psi}_{i,a}(g))+o(p^{nd})\\
    &=\sum_{a\in A}(\frac{|W|}{\phi(W)})^t\sum_{g\in G_{n-\deg(W)}^d}(\prod_{i\in [t]}\frac{\phi(W)}{|W|}\Lambda'_{W, b_i(a)}(\Tilde{\psi}_{i,a}(g))-1+1)+o(p^{nd})\\
    &=\sum_{a\in A}(\frac{|W|}{\phi(W)})^t\sum_{g\in G_{n-\deg(W)}^d}(\prod_{i\in [t]}\frac{\phi(W)}{|W|}\Lambda'_{W, b_i(a)}(\Tilde{\psi}_{i,a}(g))-1)\\
    &+\sum_{a\in A}(\frac{|W|}{\phi(W)})^t\sum_{g\in G_{n-\deg(W)}^d} 1+o(p^{nd})
\end{align*}
By using Theorem \ref{LIP}, for each $a$ we have:
\[\sum_{g\in G_{n-\deg(W)}^d}(\prod_{i\in [t]}\frac{\phi(W)}{|W|}\Lambda'_{W, b_i(a)}(\Tilde{\psi}_{i,a}(g))-1)=o_s(p^{d(n-\deg(W))})\] 
and by the definition of $\beta_W$ and Lemma \ref{beta p} :
\[\beta_W=\frac{1}{q^{\deg(W)d}}\sum_{a\in A}(\frac{|W|}{\phi(W)})^t\]
Since the Euler function is multiplicative, and by the definition of $W$, we get:
\[\beta_W=\prod_{r \text{ irr. and } \deg(r)<w(n)} \beta_r\]
Now using Lemma  \ref{beta p} we get that $\beta_W$ is bounded (with a bound that is independent of $n$),
So we have:
\begin{align*}
    &\sum_{g\in G_n^d}\prod_{i\in[t]}\Lambda'(\psi_i(g))=\sum_{a\in A}(\frac{|W|}{\phi(W)})^t\sum_{g\in G_{n-\deg(W)}^d}1+o(p^{nd})\\
    &=p^{d(n-\deg(W))}\sum_{a\in A}(\frac{|W|}{\phi(W)})^t+o(p^{nd})
\end{align*}
Now, 
\[\beta_W=\frac{1}{p^{\deg(W)d}}\sum_{a\in A}(\frac{|W|}{\phi(W)})^t\]
meaning:
\[\sum_{g\in G_n^d}\prod_{i\in[t]}\Lambda'(\psi_i(g))=p^{nd}\beta_W+o_s(p^{nd}).\]
By using the identity above we get:
\[\sum_{g\in G_n^d}\prod_{i\in[t]}\Lambda'(\psi_i(g))=p^{nd}\prod_{r \text{ irr. and } \deg(r)<w(n)} \beta_r+o_s(p^{nd})\]
and by Lemma \ref{beta p} we get:
\[\sum_{g\in G_n^d}\prod_{i\in[t]}\Lambda'(\psi_i(g))=p^{nd}\prod_{r \text{ irr. and }\deg(r)<w(n)} \beta_r+o_s(p^{nd})=p^{nd}\prod_{r \text{ irr.}} \beta_r+o_{s,t,d,L}(p^{nd})\]
And by replacing $\Lambda$ with $\Lambda'$ we can conclude Theorem \ref{LEP}.
$\square$
\newpage
\printbibliography

\end{document}